\documentclass[12pt,a4paper]{amsart}

\usepackage{amsmath,amssymb,amsfonts}
\usepackage{mathrsfs,latexsym,amsthm,stmaryrd,enumerate,epic}
\usepackage{amscd}  
\usepackage[all]{xy}

\newtheorem{theorem}{Theorem}

\newtheorem{proposition}[theorem]{Proposition}

\renewcommand{\phi}{\varphi}

\begin{document}

\title{Schur-Weyl dualities for symmetric inverse semigroups}
\author{Ganna Kudryavtseva and Volodymyr Mazorchuk}
\date{\today}

\begin{abstract}
We obtain Schur-Weyl dualities in which the algebras, acting on both 
sides, are semigroup algebras of various symmetric inverse semigroups
and their deformations.
\end{abstract}
\maketitle

\vspace{2mm}

\vspace{2mm}

\vspace{2mm}

\section{Introduction and description of results}\label{s1}

Let $V=\mathbb{C}^n$ be the natural $n$-dimensional representation of the 
group $\mathbf{GL}(n)$. Then for every $k$ the group $\mathbf{GL}(n)$ acts
diagonally on the $k$-fold tensor product $V^{\otimes k}$. At the same time
the symmetric group $S_k$ acts on $V^{\otimes k}$ by permuting the factors
of a $k$-tensor. These two actions obviously commute. Moreover, the 
classical {\em Schur-Weyl duality} from \cite{Sch1,Sch2,We} states that
$\mathbf{GL}(n)$ and $S_k$ generate full centralizers of each other on 
$V^{\otimes k}$. In particular, $\mathrm{End}_{\mathbf{GL}(n)}(V^{\otimes k})
=\mathbb{C}[S_k]$ if $n\geq k$.

There are various generalizations of the Schur-Weyl duality. In \cite{Br}
the above action of $\mathbf{GL}(n)$ was restricted to the orthogonal 
subgroup $\mathbf{O}(n)$ of $\mathbf{GL}(n)$. The corresponding centralizer 
algebra, obtained on the right hand side, is what is now known as the
{\em Brauer algebra}. Further restriction of the $\mathbf{GL}(n)$ action to 
the subgroup  $S_n$, which was considered in \cite{Jo,Ma}, gives on the right
hand side the so-called {\em partition  algebra}. Some other generalizations
are discussed in \cite{Do}, see also the references of the latter paper.

Both Brauer algebras and partition  algebras are deformations of semigroup
algebras of certain finite semigroups, which have been also intensively studied,
see for example \cite{KM,KMM,Maz} and references therein. Finite semigroups
clearly entered the game after the paper \cite{So} of L.~Solomon. In the 
latter paper it was shown that the representation $V$ of $\mathbf{GL}(n)$ can
be slightly modified such that the centralizing object obtained on the right
hand side is the {\em symmetric inverse semigroup} $\mathcal{IS}_n$, introduced
in \cite{Wa} and also known as the {\em rook monoid}, see \cite{So2}. This
idea of modification of $V$ was recently used in \cite{Gr} to obtain a
Schur-Weyl duality between $S_n$ and a generalization of the partition 
algebra, called in \cite{Gr} the {\em rook partition algebra}.

For a ``finite semigroup theorist'' there is a slight feeling of dissatisfaction
in the Schur-Weyl dualities listed above, which is explained by the fact that
the objects, appearing on the different sides of a Schur-Weyl duality, 
although closely connected to finite semigroups, still at least one 
of them has different nature. The aim of the present paper is to
establish two Schur-Weyl dualities, where on each side one has an action
of a finite inverse semigroup. 

For the first Schur-Weyl duality we have an action of the symmetric inverse
semigroup $\mathcal{IS}_n$ mentioned above on the left and an action of the 
{\em dual symmetric inverse  semigroup} $\mathcal{I}^{*}_n$ on the right.
The semigroup $\mathcal{I}^{*}_n$ was introduced in \cite{FL} as a kind of
a ``categorical dual'' for $\mathcal{IS}_n$ (see also \cite{KM2} for more 
details on the categorical approach). It is remarkable that in the present
paper the semigroup $\mathcal{I}^{*}_n$ again appears as the dual of 
$\mathcal{IS}_n$, but now with respect to a Schur-Weyl duality. 
The connection between these two types of dualities is not yet clear.
This first Schur-Weyl duality is considered in Section~\ref{s2}. 

For the second Schur-Weyl duality  we have an action of the semigroup
$\mathcal{IS}_n$ on the left and an action of the {\em partial analogue}
$\mathcal{PI}^{*}_n$ of the semigroup $\mathcal{I}^{*}_n$ on the right. The
latter semigroup was introduced in \cite{KMl} via the usual semigroup-theoretic
``partialization'' philosophy. This philosophy is rather similar to the
philosophy used to construct ``rook'' algebras. Our  second Schur-Weyl duality 
shows, in particular, an explicit connection between these philosophies. This 
second Schur-Weyl duality is considered in Section~\ref{s3}. In
Section~\ref{s4} we show that the action on the right hand side of the latter 
Schur-Weyl  duality can be ``deformed'' such that it becomes an action
of another inverse semigroup, recently constructed in \cite{Ve}.

For a semigroup $S$ we denote by $\mathbb{C}[S]$ the {\em semigroup algebra} 
of $S$ over complex numbers. If $S$ has the zero element $0$, we denote by 
$\overline{\mathbb{C}[S]}$ the {\em contracted semigroup algebra} $\mathbb{C}[S]/(0)$.
\vspace{3mm}

\noindent
{\bf Acknowledgments.} This paper was essentially written during the 
visit of the first author to Uppsala University, which was supported by 
The Royal Swedish Academy of Sciences and The Swedish Foundation for
International Cooperation in  Research and Higher Education (STINT). 
The financial support of The Academy and STINT, and the hospitality of 
Uppsala University are gratefully acknowledged. The second author is
partially supported by the Swedish Research Council.

\section{A Schur-Weyl duality for $\mathcal{IS}_n$ and
$\mathcal{I}^{*}_k$}\label{s2}

Throughout the paper $n$ and $k$ are fixed positive integers. The semigroup
$\mathcal{IS}_n$ is the semigroup of all partial injections from the set 
$\mathbf{N}=\{1,2,\dots,n\}$ to itself with respect to the usual 
composition of partial maps, see for example 
\cite[Section~2]{GM}. One can also consider the elements of $\mathcal{IS}_n$ as
bijections between different subsets of $\mathbf{N}$ (this will be important
to understand the dual nature of $\mathcal{I}^{*}_n$). A standard notation for
elements of $\mathcal{IS}_n$ and the multiplication rule in this semigroup
is best understood on the following example for $\mathcal{IS}_5$ (note that
we understand the elements of $\mathcal{IS}_n$ as maps and hence compose them
from the right to the left):
\begin{displaymath}
\left(\begin{array}{ccccc}1&2&3&4&5\\2&\varnothing&3&5&\varnothing
\end{array}\right) 
\left(\begin{array}{ccccc}1&2&3&4&5\\5&4&1&\varnothing&\varnothing
\end{array}\right) =
\left(\begin{array}{ccccc}1&2&3&4&5\\\varnothing&5&2&\varnothing&\varnothing
\end{array}\right) 
\end{displaymath}
Each element $\alpha\in \mathcal{IS}_n$ is uniquely determined by a subset
$A\subset \mathbf{N}$ and an injective map $A\to \mathbf{N}$. Abusing notation
we will denote the latter map by $\alpha$. The set $A$ is called the {\em 
domain} of $\alpha$, the set $\alpha(A)$ is called the {\em image} of $\alpha$
and the number $|A|$ is called the {\em rank} of $\alpha$, see \cite{GM} for
details. For each subset $A\subset \mathbf{N}$ we denote by $\varepsilon_A$
the idempotent of $\mathcal{IS}_n$, which corresponds to the natural
inclusion $A\hookrightarrow \mathbf{N}$. The element $\varepsilon_{\varnothing}$
is the zero element of $\mathcal{IS}_n$.

Further, the semigroup  $\mathcal{IS}_n$ can also be realized as the semigroup 
of all $n\times n$ matrices with coefficients from $\{0,1\}$ satisfying the 
condition that each row and each column of the matrix contains at most one 
non-zero entry (such matrices are called {\em rook matrices} in e.g. 
\cite{So}). The operation in the latter semigroup is the usual matrix 
multiplication. This realization defines on $V=\mathbb{C}^n$ the structure of a
$\mathbb{C}[\mathcal{IS}_n]$-module in the natural way. It is easy to see
that this module is irreducible. We call it the {\em natural representation}
of $\mathcal{IS}_n$. For each $k$ the semigroup $\mathcal{IS}_n$ acts
diagonally on the $k$-fold tensor product $V^{\otimes k}$. This is the left
hand side of our first Schur-Weyl duality.

\begin{figure}
\special{em:linewidth 0.4pt} \unitlength 0.80mm
\linethickness{0.4pt}
\begin{picture}(150.00,75.00)
\put(20.00,00.00){\makebox(0,0)[cc]{$\bullet$}}
\put(20.00,10.00){\makebox(0,0)[cc]{$\bullet$}}
\put(20.00,20.00){\makebox(0,0)[cc]{$\bullet$}}
\put(20.00,30.00){\makebox(0,0)[cc]{$\bullet$}}
\put(20.00,40.00){\makebox(0,0)[cc]{$\bullet$}}
\put(20.00,50.00){\makebox(0,0)[cc]{$\bullet$}}
\put(20.00,60.00){\makebox(0,0)[cc]{$\bullet$}}
\put(20.00,70.00){\makebox(0,0)[cc]{$\bullet$}}
\put(45.00,10.00){\makebox(0,0)[cc]{$\bullet$}}
\put(45.00,20.00){\makebox(0,0)[cc]{$\bullet$}}
\put(45.00,30.00){\makebox(0,0)[cc]{$\bullet$}}
\put(45.00,40.00){\makebox(0,0)[cc]{$\bullet$}}
\put(45.00,50.00){\makebox(0,0)[cc]{$\bullet$}}
\put(45.00,60.00){\makebox(0,0)[cc]{$\bullet$}}
\put(45.00,70.00){\makebox(0,0)[cc]{$\bullet$}}
\put(45.00,00.00){\makebox(0,0)[cc]{$\bullet$}}
\drawline(18.00,-03.00)(18.00,03.00)
\drawline(18.00,03.00)(47.00,13.00)
\drawline(47.00,13.00)(47.00,07.00)
\drawline(47.00,07.00)(18.00,-03.00)
\drawline(18.00,07.00)(18.00,23.00)
\drawline(18.00,23.00)(47.00,33.00)
\drawline(47.00,33.00)(47.00,28.00)
\drawline(47.00,28.00)(32.25,15.00)
\drawline(32.25,15.00)(47.00,03.00)
\drawline(47.00,03.00)(47.00,-03.00)
\drawline(47.00,-03.00)(18.00,07.00)
\drawline(18.00,28.00)(18.00,33.00)
\drawline(18.00,33.00)(47.00,23.00)
\drawline(47.00,23.00)(47.00,17.00)
\drawline(47.00,17.00)(18.00,28.00)
\drawline(18.00,38.00)(18.00,52.00)
\drawline(18.00,52.00)(47.00,63.00)
\drawline(47.00,63.00)(47.00,38.00)
\drawline(47.00,38.00)(18.00,38.00)
\drawline(18.00,57.00)(18.00,73.00)
\drawline(18.00,73.00)(47.00,73.00)
\drawline(47.00,73.00)(47.00,68.00)
\drawline(47.00,68.00)(18.00,57.00)
\put(70.00,00.00){\makebox(0,0)[cc]{$\bullet$}}
\put(70.00,10.00){\makebox(0,0)[cc]{$\bullet$}}
\put(70.00,20.00){\makebox(0,0)[cc]{$\bullet$}}
\put(70.00,30.00){\makebox(0,0)[cc]{$\bullet$}}
\put(70.00,40.00){\makebox(0,0)[cc]{$\bullet$}}
\put(70.00,50.00){\makebox(0,0)[cc]{$\bullet$}}
\put(70.00,60.00){\makebox(0,0)[cc]{$\bullet$}}
\put(70.00,70.00){\makebox(0,0)[cc]{$\bullet$}}
\put(95.00,10.00){\makebox(0,0)[cc]{$\bullet$}}
\put(95.00,20.00){\makebox(0,0)[cc]{$\bullet$}}
\put(95.00,30.00){\makebox(0,0)[cc]{$\bullet$}}
\put(95.00,40.00){\makebox(0,0)[cc]{$\bullet$}}
\put(95.00,50.00){\makebox(0,0)[cc]{$\bullet$}}
\put(95.00,60.00){\makebox(0,0)[cc]{$\bullet$}}
\put(95.00,70.00){\makebox(0,0)[cc]{$\bullet$}}
\put(95.00,00.00){\makebox(0,0)[cc]{$\bullet$}}
\drawline(68.00,-03.00)(68.00,03.00)
\drawline(68.00,03.00)(97.00,13.00)
\drawline(97.00,13.00)(97.00,07.00)
\drawline(97.00,07.00)(68.00,-03.00)
\drawline(68.00,07.00)(68.00,13.00)
\drawline(68.00,13.00)(97.00,23.00)
\drawline(97.00,23.00)(97.00,17.00)
\drawline(97.00,17.00)(93.00,15.00)
\drawline(93.00,15.00)(93.00,04.00)
\drawline(93.00,04.00)(97.00,03.00)
\drawline(97.00,03.00)(97.00,-03.00)
\drawline(97.00,-03.00)(68.00,07.00)
\drawline(68.00,17.00)(68.00,33.00)
\drawline(68.00,33.00)(97.00,33.00)
\drawline(97.00,33.00)(97.00,28.00)
\drawline(97.00,28.00)(68.00,17.00)
\drawline(68.00,48.00)(68.00,53.00)
\drawline(68.00,53.00)(97.00,43.00)
\drawline(97.00,43.00)(97.00,37.00)
\drawline(97.00,37.00)(68.00,48.00)
\drawline(68.00,58.00)(68.00,73.00)
\drawline(68.00,73.00)(97.00,63.00)
\drawline(97.00,63.00)(97.00,58.00)
\drawline(97.00,58.00)(68.00,58.00)
\drawline(68.00,38.00)(68.00,43.00)
\drawline(68.00,43.00)(93.00,73.00)
\drawline(93.00,73.00)(97.00,73.00)
\drawline(97.00,73.00)(97.00,68.00)
\drawline(97.00,68.00)(93.00,67.00)
\drawline(93.00,67.00)(93.00,52.00)
\drawline(93.00,52.00)(97.00,52.00)
\drawline(97.00,52.00)(97.00,47.00)
\drawline(97.00,47.00)(68.00,38.00)
\drawline(45.00,00.30)(50.00,00.30)
\drawline(55.00,00.30)(60.00,00.30)
\drawline(65.00,00.30)(70.00,00.30)
\drawline(45.00,10.30)(50.00,10.30)
\drawline(55.00,10.30)(60.00,10.30)
\drawline(65.00,10.30)(70.00,10.30)
\drawline(45.00,20.30)(50.00,20.30)
\drawline(55.00,20.30)(60.00,20.30)
\drawline(65.00,20.30)(70.00,20.30)
\drawline(45.00,30.30)(50.00,30.30)
\drawline(55.00,30.30)(60.00,30.30)
\drawline(65.00,30.30)(70.00,30.30)
\drawline(45.00,40.30)(50.00,40.30)
\drawline(55.00,40.30)(60.00,40.30)
\drawline(65.00,40.30)(70.00,40.30)
\drawline(45.00,50.30)(50.00,50.30)
\drawline(55.00,50.30)(60.00,50.30)
\drawline(65.00,50.30)(70.00,50.30)
\drawline(45.00,60.30)(50.00,60.30)
\drawline(55.00,60.30)(60.00,60.30)
\drawline(65.00,60.30)(70.00,60.30)
\drawline(45.00,70.30)(50.00,70.30)
\drawline(55.00,70.30)(60.00,70.30)
\drawline(65.00,70.30)(70.00,70.30)
\drawline(120.00,37.00)(120.00,73.00)
\drawline(120.00,73.00)(149.00,73.00)
\drawline(149.00,73.00)(149.00,37.00)
\drawline(149.00,37.00)(120.00,37.00)
\drawline(120.00,08.00)(120.00,32.00)
\drawline(120.00,32.00)(149.00,32.00)
\drawline(149.00,32.00)(149.00,27.00)
\drawline(149.00,27.00)(134.50,20.00)
\drawline(134.50,20.00)(149.00,13.00)
\drawline(149.00,13.00)(149.00,08.00)
\drawline(149.00,08.00)(120.00,08.00)
\drawline(120.00,-03.00)(120.00,03.00)
\drawline(120.00,03.00)(145.00,22.00)
\drawline(145.00,22.00)(149.00,22.00)
\drawline(149.00,22.00)(149.00,18.00)
\drawline(149.00,18.00)(132.50,03.00)
\drawline(132.50,03.00)(149.00,03.00)
\drawline(149.00,03.00)(149.00,-03.00)
\drawline(149.00,-03.00)(120.00,-03.00)
\put(122.00,00.00){\makebox(0,0)[cc]{$\bullet$}}
\put(122.00,10.00){\makebox(0,0)[cc]{$\bullet$}}
\put(122.00,20.00){\makebox(0,0)[cc]{$\bullet$}}
\put(122.00,30.00){\makebox(0,0)[cc]{$\bullet$}}
\put(122.00,40.00){\makebox(0,0)[cc]{$\bullet$}}
\put(122.00,50.00){\makebox(0,0)[cc]{$\bullet$}}
\put(122.00,60.00){\makebox(0,0)[cc]{$\bullet$}}
\put(122.00,70.00){\makebox(0,0)[cc]{$\bullet$}}
\put(147.00,10.00){\makebox(0,0)[cc]{$\bullet$}}
\put(147.00,20.00){\makebox(0,0)[cc]{$\bullet$}}
\put(147.00,30.00){\makebox(0,0)[cc]{$\bullet$}}
\put(147.00,40.00){\makebox(0,0)[cc]{$\bullet$}}
\put(147.00,50.00){\makebox(0,0)[cc]{$\bullet$}}
\put(147.00,60.00){\makebox(0,0)[cc]{$\bullet$}}
\put(147.00,70.00){\makebox(0,0)[cc]{$\bullet$}}
\put(147.00,00.00){\makebox(0,0)[cc]{$\bullet$}}
\put(57.50,35.00){\makebox(0,0)[cc]{$\cdot$}}
\put(109.50,35.00){\makebox(0,0)[cc]{$=$}}
\end{picture}
\caption{Elements of $\mathcal{I}^{*}_8$ and their
multiplication.}\label{figone}
\end{figure}
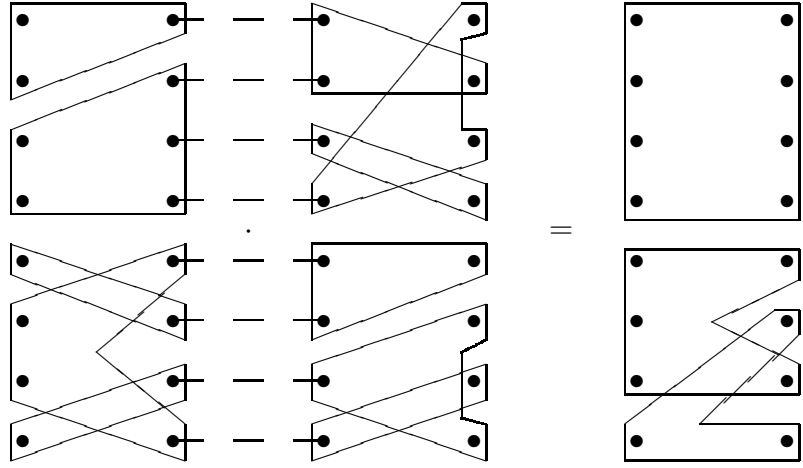

Consider the sets $\mathbf{K}=\{1,2,\dots,k\}$ and $\mathbf{K}'=
\{1',2',\dots,k'\}$. We consider ${}':\mathbf{K}\to \mathbf{K}'$ as the
{\em natural bijection} between these two sets and, abusing notation, denote
its inverse also by ${}'$ (thus $(2')'=2$). The semigroup $\mathcal{I}^{*}_k$ 
is defined in \cite{FL} as the semigroup of all bijections between different 
quotient sets of $\mathbf{K}$. Hence we can view the elements of 
$\mathcal{I}^{*}_k$ as all possible partitions of the set $\tilde{\mathbf{K}}=
\{1,2,\dots,k,1',2',\dots,k'\}$ into disjoint unions of subsets (these 
subsets will be called {\em blocks}), satisfying the condition that each block  
intersects both, $\mathbf{K}$ and $\mathbf{K}'$, non-trivially. If one drops 
the latter condition out, one obtains the list of elements of the {\em 
composition  semigroup} $\mathfrak{C}_k$, see \cite{Maz}. The multiplication on
$\mathcal{I}^{*}_k$ is much more complicated than that on $\mathcal{IS}_n$, and 
is in fact obtained by  restricting the multiplication from $\mathfrak{C}_k$. 
For an explicit formal definition of the latter we  refer the reader to 
\cite{Gr,HR,Maz}. Informally, to multiply two  partitions $\alpha$ and $\beta$ 
of $\tilde{\mathbf{K}}$ one identifies the elements $1',2',\dots,k'$ of 
$\alpha$ with the corresponding elements $1,2,\dots,k$ of $\beta$ and thus 
forms a new partition $\alpha\beta$ of $\tilde{\mathbf{K}}$ (in the latter set 
the elements  $1,2,\dots,k$ are taken from $\alpha$ and the elements 
$1',2',\dots,k'$ are taken from $\beta$), possibly deleting some ``garbage'' 
which does not contain any elements from $\mathbf{K}$ for  $\alpha$ and any 
element from $\mathbf{K}'$ for $\beta$. The partition algebra 
$\mathcal{P}_k(q)$ of \cite{Jo,Ma} is a deformation of the semigroup algebra of 
$\mathfrak{C}_k$, in which the number of garbage components, which appear
during  the above procedure, is taken into account in terms of 
a multiplicative parameter $q$. If both $\alpha$ and $\beta$ are elements 
from $\mathcal{I}^{*}_k$, then, in fact, no garbage appears. In particular, 
the algebra $\mathbb{C}[\mathcal{I}^{*}_k]$ is a subalgebra of both 
$\mathbb{C}[\mathfrak{C}_k]$ and  $\mathcal{P}_k(q)$. An example of 
multiplication of two elements from  $\mathcal{I}^{*}_8$ is given on 
Figure~\ref{figone}. Since the elements of $\mathcal{I}^{*}_k$ are defined
as certain partitions of $\tilde{\mathbf{K}}$, the set 
$\mathcal{I}^{*}_k$ is partially order with respect to inclusions in the
natural way. For $\alpha,\beta\in \mathcal{I}^{*}_k$ we will write 
$\alpha\preceq \beta$ provided that each block of the partition $\beta$
is a union of some blocks of the partition $\alpha$.

Now let us define an action of $\mathcal{I}^{*}_k$ on $V^{\otimes k}$.
Denote by $\mathbf{e}=(e_1,\dots,e_n)$ the standard basis of $V$.
For $\mathbf{i}=(i_1,\dots,i_k)\in \mathbf{N}^k$ set
\begin{displaymath}
v_{\mathbf{i}}=e_{i_1}\otimes e_{i_2}\otimes\dots\otimes e_{i_k}.
\end{displaymath}
Then the set $\mathbf{B}=\{v_{\mathbf{i}}\,:\, \mathbf{i}\in \mathbf{N}^k\}$
is a distinguished basis of $V^{\otimes k}$. For all $\alpha\in \mathfrak{C}_n$
(in particular, for all $\alpha\in\mathcal{I}^{*}_k$) and $\mathbf{i}\in
\mathbf{N}^k$ let $M(\alpha,\mathbf{i})$ denote the set of all $\mathbf{l}\in
\mathbf{N}^k$ such that for each block $\{a_1,\dots,a_p,b'_1,\dots,b'_q\}$  
of $\alpha$ (here if $\alpha\not\in \mathcal{I}^{*}_k$ it may happen that 
$p=0$ or $q=0$) we have 
\begin{displaymath}
i_{a_1}=i_{a_2}=\dots=i_{a_p}=l_{b_1}=l_{b_2}=\dots=l_{b_q}.     
\end{displaymath}
Note that $|M(\alpha,\mathbf{i})|\leq 1$ for all $\alpha\in \mathcal{I}^{*}_k$.
Now for $\alpha\in \mathfrak{C}_n$ (in particular, for all 
$\alpha\in\mathcal{I}^{*}_k$) we define that $\alpha$ acts on  $V^{\otimes k}$
as the unique linear operator such that for all $\mathbf{i}\in \mathbf{N}^k$ 
we have
\begin{equation}\label{form1}
{\alpha}(v_{\mathbf{i}})=
\begin{cases}
\displaystyle\sum_{\mathbf{l}\in M(\alpha,\mathbf{i})}  v_{\mathbf{l}},& 
M(\alpha,\mathbf{i})\neq \varnothing;\\
0,& \text{otherwise}.
\end{cases}
\end{equation}
According to \cite{Jo,Ma} this gives an action of $\mathcal{P}_k(n)$ on
$V^{\otimes k}$. By restriction, we thus also obtain an action of 
$\mathcal{I}^{*}_k$ on $V^{\otimes k}$. Now we are ready to formulate our 
first result:

\begin{theorem}\label{thm1}
\begin{enumerate}[(i)]
\item\label{thm1.1} The actions of $\mathcal{IS}_n$ and $\mathcal{I}^{*}_k$
on $V^{\otimes k}$ commute.
\item\label{thm1.2} $\mathrm{End}_{\mathcal{IS}_n}(V^{\otimes k})$ coincides
with the image of $\mathbb{C}[\mathcal{I}^{*}_k]$.
\item\label{thm1.3} $\mathrm{End}_{\mathcal{I}^*_k}(V^{\otimes k})$ 
coincides with the image of $\mathbb{C}[\mathcal{IS}_n]$.
\item\label{thm1.4} The representation of $\mathcal{IS}_n$ on 
$V^{\otimes k}$ is faithful.
\item\label{thm1.5} The representation of $\mathcal{I}^*_k$ on $V^{\otimes k}$
is faithful if and only if  $n\geq 2$ or $k=1$.
\item\label{thm1.6} The representation of
$\overline{\mathbb{C}[\mathcal{IS}_n]}$ on  $V^{\otimes k}$ is faithful 
if and only if  $k\geq n$.
\item\label{thm1.7} The representation of $\mathbb{C}[\mathcal{I}^*_k]$ on 
$V^{\otimes k}$ is faithful if and only if $k\leq n$.
\end{enumerate}
\end{theorem}

\begin{proof}
The left action of the group $S_n\subset \mathcal{IS}_n$ commutes with the
right action of $\mathcal{P}_k(n)$ on $V^{\otimes k}$ by \cite{Jo,Ma}. Hence
the left action of $S_n$ commutes with the right action of $\mathcal{I}^*_k$ 
as the latter action is obtained from the action of $\mathcal{P}_k(n)$ by
restriction. To shorten our notation, set $\varepsilon_n=
\varepsilon_{\mathbf{N}\setminus\{n\}}$. The action of $\varepsilon_n$ on 
$V$ is given by
\begin{displaymath}
\varepsilon_n(e_i)=\begin{cases}
e_i,& i\neq n;\\ 0, & i=n.\end{cases}
\end{displaymath}
Hence for any $\mathbf{i}\in\mathbf{N}^k$ we have
\begin{equation}\label{form2}
\varepsilon_n(v_{\mathbf{i}})=\begin{cases}
v_{\mathbf{i}},& n\not\in\{i_1,\dots,i_k\};\\ 
0, & \text{otherwise.}\end{cases}
\end{equation}
Now let $\alpha\in \mathcal{I}^*_k$ and $\mathbf{i}\in\mathbf{N}^k$.
Assume first that $n\not\in\{i_1,\dots,i_k\}$. As each block of the partition
$\alpha$  intersects both $\mathbf{K}$ and $\mathbf{K}'$, from the formula 
\eqref{form1} we have that $\alpha(v_{\mathbf{i}})=v_{\mathbf{j}}$, where
$\mathbf{j}\in\mathbf{N}^k$ is such that  $n\not\in \{j_1,\dots,j_k\}$.
Applying \eqref{form2} we obtain $\varepsilon_n {\alpha}(v_{\mathbf{i}})=
\alpha\varepsilon_n(v_{\mathbf{i}})$. Assume now that $n\in\{i_1,\dots,i_k\}$.
Then $\varepsilon_n(v_{\mathbf{i}})=0$ by \eqref{form2}. However, as each 
block of the partition $\alpha$ intersects both $\mathbf{K}$ and $\mathbf{K}'$,
from the formula \eqref{form1} we have that $\alpha(v_{\mathbf{i}})=
v_{\mathbf{j}}$ , where $\mathbf{j}\in\mathbf{N}^k$ is such that $n\in 
\{j_1,\dots,j_k\}$. Hence $\varepsilon_n(v_{\mathbf{j}})=0$ and we again
have $\varepsilon_n{\alpha}(v_{\mathbf{i}})=0=
\alpha\varepsilon_n(v_{\mathbf{i}})$. Therefore
$\varepsilon_n\alpha=\alpha\varepsilon_n$.  By  \cite[Theorem~3.1.4]{GM}, the
semigroup $\mathcal{IS}_n$ is  generated by its subgroup  $S_n$ and the element 
$\varepsilon_n$. The statement \eqref{thm1.1}  follows.

For $j\in\{1,2,\dots,2k\}$ let $\mathfrak{C}_k^j$ denote the set of all
elements form $\mathfrak{C}_k$, which are partitions of $\tilde{\mathbf{K}}$
into at most $j$ blocks. From  e.g. \cite[2.1]{Bl} it follows that the vector 
space $\mathrm{End}_{{S}_n}(V^{\otimes k})$ is generated by $\alpha\in
\mathfrak{C}_k^n$. Let
\begin{displaymath}
\mathbf{u}=\sum_{\alpha\in \mathfrak{C}_k^n\setminus\mathcal{I}^*_k} 
a_{\alpha}\alpha
\end{displaymath}
be a linear combination of operators acting on $V^{\otimes k}$ and 
\begin{displaymath}
 X=\{\alpha\in \mathfrak{C}_k^n\setminus\mathcal{I}^*_k:a_{\alpha}\neq 0\}.
\end{displaymath}
We would like to show that the condition $\mathbf{u}\varepsilon_n=
\varepsilon_n\mathbf{u}$ implies $X=\varnothing$. Assume $X\neq\varnothing$.
Let $\alpha\in X$ be a minimal element with respect to $\preceq$, that is
a partition, which does not properly contain any other partition from $X$. As
$\alpha\not\in \mathcal{I}^*_k$, the element $\alpha$ contains some block, 
say $B$, which is contained in either $\mathbf{K}$ or $\mathbf{K}'$. Consider
some  map $f:\tilde{\mathbf{K}}\to \mathbf{N}$, which satisfies the following
conditions:
\begin{itemize} 
\item $f$ is constant on blocks of $\alpha$;
\item $f$ has different values on elements from different blocks;
\item $f$ has  value $n$ on elements from the block $B$.
\end{itemize}
Such map exists because $\alpha\in \mathfrak{C}_k^n$. Consider now the elements
\begin{equation}\label{elements}
{v}=e_{f(1)}\otimes e_{f(2)}\otimes\dots\otimes e_{f(k)}\quad\text{ and }\quad
{w}=e_{f(1')}\otimes e_{f(2')}\otimes\dots\otimes e_{f(k')}.
\end{equation}
Assume first that $B\subset \mathbf{K}$. Then $\varepsilon_n(v)=0$ as
$n$ occurs among $f(1),\dots,f(k)$ and hence $\mathbf{u}\varepsilon_n(v)=0$
as well. On the other hand, the element $\alpha(v)$, when expressed as a linear combination of elements from $\mathbf{B}$, has a non-zero coefficient at $w$
because of \eqref{form1} and the definition of $f$. Further, for any 
$\alpha'\in X$ different from $\alpha$ the coefficient of $\alpha'(v)$ at $w$ 
(again when $\alpha'(v)$ is expressed as a linear combination of the elements 
from $\mathbf{B}$) is zero because of \eqref{form1}, the minimality of the 
partition $\alpha$ with respect to $\preceq$ and the definition of $f$. Hence 
the element $\mathbf{u}(v)$, when expressed as a linear combination of the elements from $\mathbf{B}$, has a non-zero coefficient at $w$. But $n$ does not 
occur among $f(1'),\dots,f(k')$, which means $\varepsilon_n(w)=w$. As the 
action of $\varepsilon_n$ is diagonal  with respect to the basis $\mathbf{B}$, 
we get that $\varepsilon_n\mathbf{u}(v) \neq 0$. 

In the case $B\subset \mathbf{N}'$ by similar arguments we get 
$\mathbf{u}\varepsilon_n(v)\neq 0$ while $\varepsilon_n\mathbf{u}(v)=0$. Hence
$\mathbf{u}\varepsilon_n\neq \varepsilon_n\mathbf{u}$, which means that 
$\mathrm{End}_{\mathcal{IS}_n}(V^{\otimes k})$ is already generated by 
$\mathfrak{C}_k^n\cap\mathcal{I}^*_k$. The statement \eqref{thm1.2} follows.

The statement \eqref{thm1.3} is now a standard double-centralizer property. 
As the semigroup $\mathcal{IS}_n$ is an inverse semigroup
(see e.g. \cite[Theorem~2.6.7]{GM}), the semigroup algebra
$\mathbb{C}[\mathcal{IS}_n]$ is semisimple by \cite[Theorem~4.4]{Mu}.
Hence the image of $\mathbb{C}[\mathcal{IS}_n]$ in
$\mathrm{End}_{\mathbb{C}}(V^{\otimes k})$ is semisimple as well. The
double-centralizer property for semisimple algebras is a trivial case
of Tachikawa's theory of dominance dimension, see \cite{Ta}
(the idea of the above argument is taken from \cite[Theorem~2.3]{Bl} 
and \cite[Theorem~2.8]{KSX}). The statement \eqref{thm1.3} follows. 

The representation of $\mathcal{IS}_n$ on $V$ is faithful by the definition. 
Hence for any $\pi,\tau\in\mathcal{IS}_n$, $\pi\neq\tau$, there exists 
$i\in\mathbf{N}$ such that  $\pi e_i\neq \tau e_i$. But then 
$\pi(e_i\otimes\dots\otimes e_i)\neq\tau(e_i\otimes\dots\otimes e_i)$ and hence
the actions of $\pi$ and $\tau$ on $V^{\otimes k}$ are different. This proves 
\eqref{thm1.4}.

As $|\mathcal{I}^*_1|=1$, any representation of $\mathcal{I}^*_1$ 
is faithful. If $n=1$ and $k>1$, then $|\mathcal{I}^*_k|>1$. However, the 
formula \eqref{form1} says that all  elements of $\mathcal{I}^*_k$ are 
represented by the identity operator on $V^{\otimes k}$. Hence the
representation of $\mathcal{I}^*_k$ on $V^{\otimes k}$ is not faithful in the
case $n=1$ and $k>1$. Finally, assume that $n>1$. Let $\alpha,\beta\in 
\mathcal{I}^*_k$ and assume that the actions of $\alpha$ and $\beta$ on 
$V^{\otimes k}$ coincide. Let $B$ be some block of $\alpha$. Consider the map 
$f:\tilde{\mathbf{K}}\to  \mathbf{N}$ defined as follows:
\begin{displaymath}
f(x)=\begin{cases}1,&x\in B;\\2,&x\not\in B,\end{cases}
\end{displaymath}
and consider the corresponding elements $v$ and $w$ given by \eqref{elements}.
By \eqref{form1} and our construction of $f$, we have $\alpha(v)=w$. As
$\alpha(v)=\beta(v)$, from \eqref{form1} it follows that $B$ should be a 
union of blocks of $\beta$. Hence each block of $\alpha$ is a union of some
blocks of $\beta$. Analogously, each block of $\beta$ is a union of some 
blocks of $\alpha$. This implies $\alpha=\beta$ and proves \eqref{thm1.5}.

To prove \eqref{thm1.6} we first note that the zero element 
$\varepsilon_{\varnothing}$ of $\mathcal{IS}_n$ acts as the zero operator
on $V^{\otimes k}$. Hence $V^{\otimes k}$ is even an
$\overline{\mathbb{C}[\mathcal{IS}_n]}$-module. If $k<n$ then from e.g.
\cite[Theorem~3.22(a)]{HR} it follows that already the restriction of 
the $\overline{\mathbb{C}[\mathcal{IS}_n]}$-action to
$\mathbb{C}[S_n]$ does not give a faithful representation of 
$\mathbb{C}[S_n]$, as not all simple $\mathbb{C}[S_n]$-modules occur in the
decomposition of $V^{\otimes k}$. Hence for $k<n$ the action of 
$\overline{\mathbb{C}[\mathcal{IS}_n]}$ on $V^{\otimes k}$ is not faithful.

Let now $k\geq n$. Consider some linear combination
\begin{displaymath}
\mathbf{u}=\sum_{\alpha\in\mathcal{IS}_n\setminus\{\varepsilon_{\varnothing}\}}
a_{\alpha}\alpha 
\end{displaymath}
and assume that $\mathbf{u}$ annihilates $V^{\otimes k}$. Then, in particular, 
$\mathbf{u}v=0$, where 
\begin{displaymath}
v=e_1\otimes e_2\otimes\dots\otimes e_{n-1}\otimes e_n
\otimes e_n\otimes \dots \otimes e_n. 
\end{displaymath}
However, the element $v$ is annihilated by all elements
$\alpha\in \mathcal{IS}_n$ of rank at most $n-1$. At the same time
the elemens of  rank $n$ map $v$ to linearly independent elements of
$V^{\otimes k}$. This implies that $a_{\alpha}=0$ for all $\alpha$ of rank $n$.
Applying now exactly the same arguments to the vector 
\begin{displaymath}
v'=e_1\otimes e_2\otimes\dots\otimes e_{n-2}\otimes e_{n-1}\otimes
e_{n-1}\otimes \dots  \otimes e_{n-1} 
\end{displaymath}
we obtain that $a_{\alpha}=0$ for all $\alpha$ with domain $\{1,2,\dots,n-1\}$.
Analogously one shows that $a_{\alpha}=0$ for all $\alpha$ of rank $n-1$.
Proceeding by induction on the rank of $\alpha$ we thus get $a_{\alpha}=0$ 
for all $\alpha$. This proves the statement \eqref{thm1.6}.

Finally, let us prove the statement \eqref{thm1.7}. If $k>n$, then we recall 
that during the proof of \eqref{thm1.2} we saw that the image of 
$\mathbb{C}[\mathcal{I}^*_k]$ in $\mathrm{End}_{\mathbb{C}}(V^{\otimes k})$ is 
generated already by the image of  $\mathfrak{C}_k^n\cap\mathcal{I}^*_k$. Hence 
for $k>n$ the representation of $\mathbb{C}[\mathcal{I}^*_k]$ on 
$V^{\otimes k}$ is not faithful. 

Let $k\leq n$. Consider some linear combination
\begin{displaymath}
\mathbf{u}=\sum_{\alpha\in\mathcal{I}^*_k} a_{\alpha}\alpha 
\end{displaymath}
and assume that $\mathbf{u}$ annihilates $V^{\otimes k}$. Let  $\alpha\in
\mathcal{I}^*_k$ be minimal with respect to the partial order $\preceq$.
Consider some  map $f:\tilde{\mathbf{K}}\to \mathbf{N}$, which satisfies the
following conditions:
\begin{itemize}
\item $f$ is constant on  blocks  of $\alpha$;
\item $f$ has different values on elements from  different blocks.
\end{itemize}
Such a map exists as $\alpha$ has at most $k$ blocks and $k\leq n$. 
Consider now the corresponding elements $v$ and $w$ given by \eqref{elements}.
From \eqref{form1} we get that  $\alpha(v)=w$ while the coefficient of 
$\beta(v)$ at $w$ (when expressed with respect to $\mathbf{B}$) is zero for 
all $\beta\in \mathcal{I}^*_k$, $\beta\neq \alpha$, because of the minimality
of $\alpha$ and the choice of $f$. Since $\mathbf{u}(v)=0$, we thus must have
$a_{\alpha}=0$. Proceeding in the same way with respect to the partial order 
$\preceq$ on $\mathcal{I}^*_k$ we obtain $a_{\alpha}=0$ for all $\alpha\in 
\mathcal{I}^*_k$ and the statement \eqref{thm1.7} follows. This completes the 
proof of the theorem.
\end{proof}

\section{A Schur-Weyl duality for $\mathcal{IS}_n$ and
$\mathcal{PI}^{*}_k$}\label{s3}

For the second Schur-Weyl duality we consider the trivial 
$\mathcal{IS}_n$-module $\mathbb{C}$ on which all elements of $\mathcal{IS}_n$
(including the zero element $\varepsilon_{\varnothing}$) act via the identity
transformation. We denote by $e_0$ some basis element of $\mathbb{C}$.
Consider now the $\mathcal{IS}_n$-module $U=V\oplus \mathbb{C}$ and the 
vector space  $U^{\otimes k}$ as an $\mathcal{IS}_n$-module with respect to 
the diagonal action of  $\mathcal{IS}_n$. This is the left hand side.

\begin{figure}
\special{em:linewidth 0.4pt} \unitlength 0.80mm
\linethickness{0.4pt}
\begin{picture}(150.00,75.00)
\put(20.00,00.00){\makebox(0,0)[cc]{$\bullet$}}
\put(20.00,10.00){\makebox(0,0)[cc]{$\bullet$}}
\put(20.00,20.00){\makebox(0,0)[cc]{$\bullet$}}
\put(20.00,30.00){\makebox(0,0)[cc]{$\bullet$}}
\put(20.00,40.00){\makebox(0,0)[cc]{$\bullet$}}
\put(20.00,50.00){\makebox(0,0)[cc]{$\bullet$}}
\put(20.00,60.00){\makebox(0,0)[cc]{$\bullet$}}
\put(20.00,70.00){\makebox(0,0)[cc]{$\bullet$}}
\put(45.00,10.00){\makebox(0,0)[cc]{$\bullet$}}
\put(45.00,20.00){\makebox(0,0)[cc]{$\bullet$}}
\put(45.00,30.00){\makebox(0,0)[cc]{$\bullet$}}
\put(45.00,40.00){\makebox(0,0)[cc]{$\bullet$}}
\put(45.00,50.00){\makebox(0,0)[cc]{$\bullet$}}
\put(45.00,60.00){\makebox(0,0)[cc]{$\bullet$}}
\put(45.00,70.00){\makebox(0,0)[cc]{$\bullet$}}
\put(45.00,00.00){\makebox(0,0)[cc]{$\bullet$}}
\drawline(18.00,-03.00)(18.00,13.00)
\drawline(18.00,13.00)(22.00,13.00)
\drawline(22.00,13.00)(43.00,23.00)
\drawline(43.00,23.00)(47.00,23.00)
\drawline(47.00,23.00)(47.00,-03.00)
\drawline(47.00,-03.00)(18.00,-03.00)
\drawline(18.00,18.00)(18.00,32.00)
\drawline(18.00,32.00)(47.00,32.00)
\drawline(47.00,32.00)(47.00,27.00)
\drawline(47.00,27.00)(22.00,18.00)
\drawline(22.00,18.00)(18.00,18.00)
\drawline(18.00,38.00)(18.00,42.00)
\drawline(18.00,42.00)(43.00,52.00)
\drawline(43.00,52.00)(47.00,52.00)
\drawline(47.00,52.00)(47.00,48.00)
\drawline(47.00,48.00)(22.00,38.00)
\drawline(22.00,38.00)(18.00,38.00)
\drawline(18.00,48.00)(18.00,62.00)
\drawline(18.00,62.00)(43.00,72.00)
\drawline(43.00,72.00)(47.00,72.00)
\drawline(47.00,72.00)(47.00,58.00)
\drawline(47.00,58.00)(22.00,48.00)
\drawline(22.00,48.00)(18.00,48.00)
\drawline(18.00,68.00)(18.00,72.00)
\drawline(18.00,72.00)(22.00,72.00)
\drawline(22.00,72.00)(22.00,68.00)
\drawline(22.00,68.00)(18.00,68.00)
\drawline(43.00,38.00)(43.00,42.00)
\drawline(43.00,42.00)(47.00,42.00)
\drawline(47.00,42.00)(47.00,38.00)
\drawline(47.00,38.00)(43.00,38.00)
\drawline(68.00,-03.00)(68.00,03.00)
\drawline(68.00,03.00)(91.00,22.00)
\drawline(91.00,22.00)(97.00,22.00)
\drawline(97.00,22.00)(97.00,18.00)
\drawline(97.00,18.00)(91.00,18.00)
\drawline(91.00,18.00)(91.00,03.00)
\drawline(91.00,03.00)(97.00,03.00)
\drawline(97.00,03.00)(97.00,-03.00)
\drawline(97.00,-03.00)(68.00,-03.00)
\drawline(68.00,08.00)(68.00,23.00)
\drawline(68.00,23.00)(93.00,32.00)
\drawline(93.00,32.00)(97.00,32.00)
\drawline(97.00,32.00)(97.00,28.00)
\drawline(97.00,28.00)(93.00,28.00)
\drawline(93.00,28.00)(72.00,08.00)
\drawline(72.00,08.00)(68.00,08.00)
\drawline(68.00,28.00)(68.00,32.00)
\drawline(68.00,32.00)(74.00,32.00)
\drawline(74.00,32.00)(74.00,48.00)
\drawline(74.00,48.00)(68.00,48.00)
\drawline(68.00,48.00)(68.00,52.00)
\drawline(68.00,52.00)(97.00,52.00)
\drawline(97.00,52.00)(97.00,38.00)
\drawline(97.00,38.00)(93.00,38.00)
\drawline(93.00,38.00)(72.00,28.00)
\drawline(72.00,28.00)(68.00,28.00)
\drawline(68.00,68.00)(68.00,72.00)
\drawline(68.00,72.00)(97.00,72.00)
\drawline(97.00,72.00)(97.00,58.00)
\drawline(97.00,58.00)(93.00,58.00)
\drawline(93.00,58.00)(68.00,68.00)
\drawline(68.00,38.00)(68.00,42.00)
\drawline(68.00,42.00)(72.00,42.00)
\drawline(72.00,42.00)(72.00,38.00)
\drawline(72.00,38.00)(68.00,38.00)
\drawline(68.00,58.00)(68.00,62.00)
\drawline(68.00,62.00)(72.00,62.00)
\drawline(72.00,62.00)(72.00,58.00)
\drawline(72.00,58.00)(68.00,58.00)
\drawline(120.00,-03.00)(120.00,12.00)
\drawline(120.00,12.00)(124.00,12.00)
\drawline(124.00,12.00)(145.00,32.00)
\drawline(145.00,32.00)(149.00,32.00)
\drawline(149.00,32.00)(149.00,18.00)
\drawline(149.00,18.00)(143.00,18.00)
\drawline(143.00,18.00)(143.00,02.00)
\drawline(143.00,02.00)(149.00,02.00)
\drawline(149.00,02.00)(149.00,-03.00)
\drawline(149.00,-03.00)(120.00,-03.00)
\drawline(120.00,18.00)(120.00,42.00)
\drawline(120.00,42.00)(145.00,52.00)
\drawline(145.00,52.00)(149.00,52.00)
\drawline(149.00,52.00)(149.00,38.00)
\drawline(149.00,38.00)(145.00,38.00)
\drawline(145.00,38.00)(124.00,18.00)
\drawline(124.00,18.00)(120.00,18.00)
\drawline(120.00,48.00)(120.00,52.00)
\drawline(120.00,52.00)(124.00,52.00)
\drawline(124.00,52.00)(124.00,48.00)
\drawline(124.00,48.00)(120.00,48.00)
\drawline(120.00,58.00)(120.00,62.00)
\drawline(120.00,62.00)(124.00,62.00)
\drawline(124.00,62.00)(124.00,58.00)
\drawline(124.00,58.00)(120.00,58.00)
\drawline(120.00,68.00)(120.00,72.00)
\drawline(120.00,72.00)(124.00,72.00)
\drawline(124.00,72.00)(124.00,68.00)
\drawline(124.00,68.00)(120.00,68.00)
\drawline(145.00,08.00)(145.00,12.00)
\drawline(145.00,12.00)(149.00,12.00)
\drawline(149.00,12.00)(149.00,08.00)
\drawline(149.00,08.00)(145.00,08.00)
\drawline(145.00,58.00)(145.00,62.00)
\drawline(145.00,62.00)(149.00,62.00)
\drawline(149.00,62.00)(149.00,58.00)
\drawline(149.00,58.00)(145.00,58.00)
\drawline(145.00,68.00)(145.00,72.00)
\drawline(145.00,72.00)(149.00,72.00)
\drawline(149.00,72.00)(149.00,68.00)
\drawline(149.00,68.00)(145.00,68.00)
\drawline(93.00,08.00)(93.00,12.00)
\drawline(93.00,12.00)(97.00,12.00)
\drawline(97.00,12.00)(97.00,08.00)
\drawline(97.00,08.00)(93.00,08.00)
\put(70.00,00.00){\makebox(0,0)[cc]{$\bullet$}}
\put(70.00,10.00){\makebox(0,0)[cc]{$\bullet$}}
\put(70.00,20.00){\makebox(0,0)[cc]{$\bullet$}}
\put(70.00,30.00){\makebox(0,0)[cc]{$\bullet$}}
\put(70.00,40.00){\makebox(0,0)[cc]{$\bullet$}}
\put(70.00,50.00){\makebox(0,0)[cc]{$\bullet$}}
\put(70.00,60.00){\makebox(0,0)[cc]{$\bullet$}}
\put(70.00,70.00){\makebox(0,0)[cc]{$\bullet$}}
\put(95.00,10.00){\makebox(0,0)[cc]{$\bullet$}}
\put(95.00,20.00){\makebox(0,0)[cc]{$\bullet$}}
\put(95.00,30.00){\makebox(0,0)[cc]{$\bullet$}}
\put(95.00,40.00){\makebox(0,0)[cc]{$\bullet$}}
\put(95.00,50.00){\makebox(0,0)[cc]{$\bullet$}}
\put(95.00,60.00){\makebox(0,0)[cc]{$\bullet$}}
\put(95.00,70.00){\makebox(0,0)[cc]{$\bullet$}}
\put(95.00,00.00){\makebox(0,0)[cc]{$\bullet$}}
\put(122.00,00.00){\makebox(0,0)[cc]{$\bullet$}}
\put(122.00,10.00){\makebox(0,0)[cc]{$\bullet$}}
\put(122.00,20.00){\makebox(0,0)[cc]{$\bullet$}}
\put(122.00,30.00){\makebox(0,0)[cc]{$\bullet$}}
\put(122.00,40.00){\makebox(0,0)[cc]{$\bullet$}}
\put(122.00,50.00){\makebox(0,0)[cc]{$\bullet$}}
\put(122.00,60.00){\makebox(0,0)[cc]{$\bullet$}}
\put(122.00,70.00){\makebox(0,0)[cc]{$\bullet$}}
\put(147.00,10.00){\makebox(0,0)[cc]{$\bullet$}}
\put(147.00,20.00){\makebox(0,0)[cc]{$\bullet$}}
\put(147.00,30.00){\makebox(0,0)[cc]{$\bullet$}}
\put(147.00,40.00){\makebox(0,0)[cc]{$\bullet$}}
\put(147.00,50.00){\makebox(0,0)[cc]{$\bullet$}}
\put(147.00,60.00){\makebox(0,0)[cc]{$\bullet$}}
\put(147.00,70.00){\makebox(0,0)[cc]{$\bullet$}}
\put(147.00,00.00){\makebox(0,0)[cc]{$\bullet$}}
\drawline(45.00,00.30)(50.00,00.30)
\drawline(55.00,00.30)(60.00,00.30)
\drawline(65.00,00.30)(70.00,00.30)
\drawline(45.00,10.30)(50.00,10.30)
\drawline(55.00,10.30)(60.00,10.30)
\drawline(65.00,10.30)(70.00,10.30)
\drawline(45.00,20.30)(50.00,20.30)
\drawline(55.00,20.30)(60.00,20.30)
\drawline(65.00,20.30)(70.00,20.30)
\drawline(45.00,30.30)(50.00,30.30)
\drawline(55.00,30.30)(60.00,30.30)
\drawline(65.00,30.30)(70.00,30.30)
\drawline(45.00,40.30)(50.00,40.30)
\drawline(55.00,40.30)(60.00,40.30)
\drawline(65.00,40.30)(70.00,40.30)
\drawline(45.00,50.30)(50.00,50.30)
\drawline(55.00,50.30)(60.00,50.30)
\drawline(65.00,50.30)(70.00,50.30)
\drawline(45.00,60.30)(50.00,60.30)
\drawline(55.00,60.30)(60.00,60.30)
\drawline(65.00,60.30)(70.00,60.30)
\drawline(45.00,70.30)(50.00,70.30)
\drawline(55.00,70.30)(60.00,70.30)
\drawline(65.00,70.30)(70.00,70.30)
\put(57.50,35.00){\makebox(0,0)[cc]{$\cdot$}}
\put(109.50,35.00){\makebox(0,0)[cc]{$=$}}
\end{picture}
\caption{Elements of $\mathcal{PI}^*_8$ and their
multiplication.}\label{figtwo}
\end{figure}
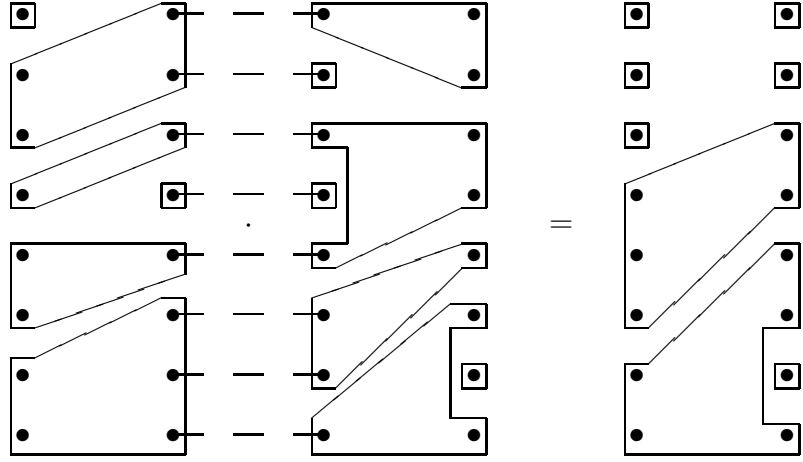

To describe the right hand side we have to define another semigroup, namely 
the {\em partial dual inverse symmetric semigroup} $\mathcal{PI}^*_k$. This
semigroup was introduced in \cite{KMl}. The elements of $\mathcal{PI}^*_k$
are all possible partitions $\alpha$ of subsets 
$A\subset \tilde{\mathbf{K}}$, which satisfy the condition that each block 
of $\alpha$ has a non-trivial intersection with both $\mathbf{K}$ and
$\mathbf{K}'$. We can consider  $\mathcal{PI}^*_k$ as a subset of 
$\mathfrak{C}_k$ extending each $\alpha\in \mathcal{PI}^*_k$ to a
partition of $\tilde{\mathbf{K}}$ as follows: if $\alpha$ is a partition of
some $A\subset \tilde{\mathbf{K}}$, then we add to this partition all elements
from $\tilde{\mathbf{K}}\setminus A$ as separate one-element blocks. In this
way $\mathcal{PI}^*_k$ becomes a subset, but not a subsemigroup of
$\mathfrak{C}_k$ (an example, illustrating that  $\mathcal{PI}^*_k$ is not 
closed with respect to the multiplication on $\mathfrak{C}_k$, can be found on
\cite[Figure~2]{KMl}). To make $\mathcal{PI}^*_k$ into a subsemigroup the
multiplication should be changed as follows: Let $\alpha,\beta\in
\mathcal{PI}^*_k$. Identify $\mathbf{K}'$-elements of $\alpha$ with the
corresponding $\mathbf{K}$-elements of $\beta$. Now those blocks, which 
do not contain any one-element blocks from $\alpha$ or $\beta$ survive, and
all other blocks break down into one-element blocks. We refer the reader
to \cite[2.1]{KMl} for the formal definition. An example of multiplication
of two elements from $\mathcal{PI}^*_8$ is given on Figure~\ref{figtwo}. The
algebra $\mathbb{C}[\mathcal{PI}^*_k]$ is a subalgebra of the rook partition 
algebra from \cite[2.1]{Gr} in the natural way. This follows immediately by 
comparing the definitions.

Now let us define an action of $\mathcal{PI}^*_k$ on $U^{\otimes k}$. This
follows closely \cite[Section~5]{So} and \cite[Section~2]{Gr}. The set
$\mathbf{B}'=\{v_{\mathbf{i}}\,:\,\mathbf{i}\in (\mathbf{N}\cup\{0\})^k\}$ is 
a distinguished basis of $U^{\otimes k}$. For $\alpha\in \mathcal{PI}^{*}_k$ 
and $\mathbf{i}\in (\mathbf{N}\cup\{0\})^k$ let $M(\alpha,\mathbf{i})$ denote 
the set of all $\mathbf{l}\in (\mathbf{N}\cup\{0\})^k$ such
that the following two conditions are satisfied:
\begin{itemize}
\item for each block $\{a_1,\dots,a_p,b'_1,\dots,b'_q\}$ of $\alpha$ we  have 
\begin{displaymath}
i_{a_1}=i_{a_2}=\dots=i_{a_p}=l_{b_1}=l_{b_2}=\dots=l_{b_q};    
\end{displaymath}
\item for any $a\in \mathbf{K}$ and $b\in \mathbf{K}'$ which do not belong 
to any block of $\alpha$ we have $i_a=0=l_b$.
\end{itemize}
Again note that $|M(\alpha,\mathbf{i})|\leq 1$ for all $\alpha\in 
\mathcal{PI}^{*}_k$ and $\mathbf{i}\in (\mathbf{N}\cup\{0\})^k$.
Now for $\alpha\in \mathcal{PI}^{*}_k$ we define the action of  $\alpha$ on  
$V^{\otimes k}$  via the formula \eqref{form1}. This action is in fact the 
restriction of the action of the rook partition algebra, constructed in 
\cite[2.1]{Gr}. In particular, we automatically obtain an action of the
semigroup $\mathcal{PI}^{*}_k$ on $U^{\otimes k}$.  Now we are  ready to
formulate our next  result.

\begin{theorem}\label{thm2}
\begin{enumerate}[(i)]
\item\label{thm2.1} The actions of $\mathcal{IS}_n$ and $\mathcal{PI}^{*}_k$
on $U^{\otimes k}$ commute.
\item\label{thm2.2} $\mathrm{End}_{\mathcal{IS}_n}(U^{\otimes k})$ coincides
with the image of $\mathbb{C}[\mathcal{PI}^{*}_k]$.
\item\label{thm2.3} $\mathrm{End}_{\mathcal{PI}^*_k}(U^{\otimes k})$
coincides with the image of $\mathbb{C}[\mathcal{IS}_n]$.
\item\label{thm2.4} The representation of $\mathcal{IS}_n$ on 
$U^{\otimes k}$ is faithful.
\item\label{thm2.5} The representation of $\mathcal{PI}^*_k$ on 
$U^{\otimes k}$ is faithful.
\item\label{thm2.6} The representation of $\mathbb{C}[\mathcal{IS}_n]$ 
on  $U^{\otimes k}$ is faithful if and only if  $k\geq n$.
\item\label{thm2.7} The representation of $\mathbb{C}[\mathcal{PI}^*_k]$ on 
$U^{\otimes k}$ is faithful if and only if $k\leq n$.
\end{enumerate}
\end{theorem}

\begin{proof}
By \cite[Theorem~18]{Gr}, the right action of $\mathcal{PI}^*_k$ on 
$U^{\otimes k}$ commutes with the left action of $S_n$. So, it is enough 
to check that the right action of $\mathcal{PI}^*_k$ commutes with the
left action of $\varepsilon_n$. This is a straightforward calculation using
definitions, which is similar to one in the proof of 
Theorem~\ref{thm1}\eqref{thm1.1}. This proves the statement \eqref{thm2.1}.

Analogously, because of \cite[Theorem~18]{Gr} the proof of the statement
\eqref{thm2.2} is similar to that of Theorem~\ref{thm1}\eqref{thm1.2}.
The proof of the statement \eqref{thm2.3} is exactly the same as one of 
Theorem~\ref{thm1}\eqref{thm1.3}. As $V^{\otimes k}$ is a submodule of
$U^{\otimes k}$, the statement \eqref{thm2.4} follows from 
Theorem~\ref{thm1}\eqref{thm1.4}. The statement \eqref{thm2.5}
is proved analogously to Theorem~\ref{thm1}\eqref{thm1.5}.

As $\mathbb{C}$ is the trivial $\mathcal{IS}_n$-module, by the additivity 
of the tensor product the module $U^{\otimes k}$ decomposes into a direct 
sum of $\mathcal{IS}_n$-modules, each of which is isomorphic to some
$V^{\otimes r}$, where $r\leq k$. Hence the fact that the representation of
$\mathbb{C}[\mathcal{IS}_n]$  on  $U^{\otimes k}$ is not faithful for 
$k< n$ follows from Theorem~\ref{thm1}\eqref{thm1.6}. If $k\geq n$, then
$V^{\otimes k}$ is a direct summand of $U^{\otimes k}$ and hence the
representation of $\overline{\mathbb{C}[\mathcal{IS}_n]}$  on  
$U^{\otimes k}$ is faithful by Theorem~\ref{thm1}\eqref{thm1.6}. At the
same time the action of $\varepsilon_{\varnothing}$ on $U^{\otimes k}$ is 
obviously non-zero as the trivial $\mathcal{IS}_n$-module is a direct summand
of $U^{\otimes k}$ as well. This implies \eqref{thm2.6}.

Finally, for $k>n$ the fact that the representation of
$\mathbb{C}[\mathcal{PI}^*_k]$ on  $U^{\otimes k}$ is not faithful follows
from \cite[Theorem~9]{Gr}. For $k\leq n$ the fact that the representation of
$\mathbb{C}[\mathcal{PI}^*_k]$ on  $U^{\otimes k}$ is faithful is proved
analogously to the corresponding part of Theorem~\ref{thm1}\eqref{thm1.7}.
\end{proof}

\section{Deformations of the second Schur-Weyl duality}\label{s4}

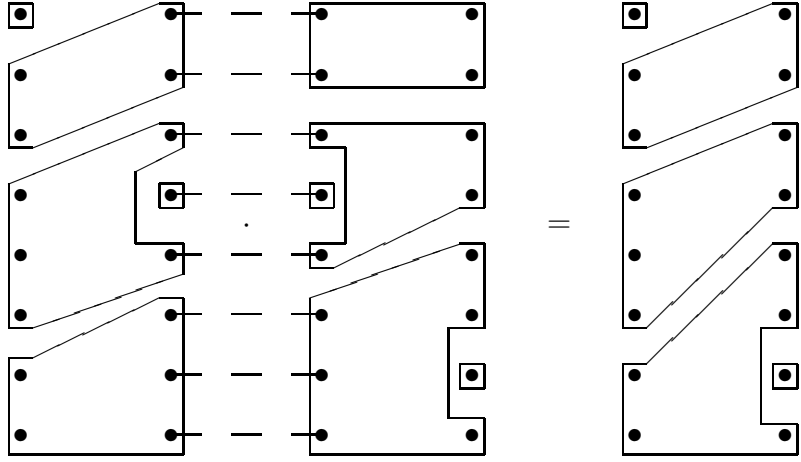
\begin{figure}
\special{em:linewidth 0.4pt} \unitlength 0.80mm
\linethickness{0.4pt}
\begin{picture}(150.00,75.00)
\put(20.00,00.00){\makebox(0,0)[cc]{$\bullet$}}
\put(20.00,10.00){\makebox(0,0)[cc]{$\bullet$}}
\put(20.00,20.00){\makebox(0,0)[cc]{$\bullet$}}
\put(20.00,30.00){\makebox(0,0)[cc]{$\bullet$}}
\put(20.00,40.00){\makebox(0,0)[cc]{$\bullet$}}
\put(20.00,50.00){\makebox(0,0)[cc]{$\bullet$}}
\put(20.00,60.00){\makebox(0,0)[cc]{$\bullet$}}
\put(20.00,70.00){\makebox(0,0)[cc]{$\bullet$}}
\put(45.00,10.00){\makebox(0,0)[cc]{$\bullet$}}
\put(45.00,20.00){\makebox(0,0)[cc]{$\bullet$}}
\put(45.00,30.00){\makebox(0,0)[cc]{$\bullet$}}
\put(45.00,40.00){\makebox(0,0)[cc]{$\bullet$}}
\put(45.00,50.00){\makebox(0,0)[cc]{$\bullet$}}
\put(45.00,60.00){\makebox(0,0)[cc]{$\bullet$}}
\put(45.00,70.00){\makebox(0,0)[cc]{$\bullet$}}
\put(45.00,00.00){\makebox(0,0)[cc]{$\bullet$}}
\drawline(18.00,-03.00)(18.00,13.00)
\drawline(18.00,13.00)(22.00,13.00)
\drawline(22.00,13.00)(43.00,23.00)
\drawline(43.00,23.00)(47.00,23.00)
\drawline(47.00,23.00)(47.00,-03.00)
\drawline(47.00,-03.00)(18.00,-03.00)
\drawline(18.00,18.00)(18.00,32.00)
\drawline(39.00,32.00)(47.00,32.00)
\drawline(47.00,32.00)(47.00,27.00)
\drawline(47.00,27.00)(22.00,18.00)
\drawline(22.00,18.00)(18.00,18.00)
\drawline(18.00,18.00)(18.00,42.00)
\drawline(18.00,42.00)(43.00,52.00)
\drawline(43.00,52.00)(47.00,52.00)
\drawline(47.00,52.00)(47.00,48.00)
\drawline(47.00,48.00)(39.00,44.00)
\drawline(39.00,44.00)(39.00,32.00)
\drawline(18.00,48.00)(18.00,62.00)
\drawline(18.00,62.00)(43.00,72.00)
\drawline(43.00,72.00)(47.00,72.00)
\drawline(47.00,72.00)(47.00,58.00)
\drawline(47.00,58.00)(22.00,48.00)
\drawline(22.00,48.00)(18.00,48.00)
\drawline(18.00,68.00)(18.00,72.00)
\drawline(18.00,72.00)(22.00,72.00)
\drawline(22.00,72.00)(22.00,68.00)
\drawline(22.00,68.00)(18.00,68.00)
\drawline(43.00,38.00)(43.00,42.00)
\drawline(43.00,42.00)(47.00,42.00)
\drawline(47.00,42.00)(47.00,38.00)
\drawline(47.00,38.00)(43.00,38.00)
\drawline(68.00,-03.00)(68.00,03.00)
\drawline(68.00,03.00)(68.00,18.00)
\drawline(97.00,22.00)(97.00,18.00)
\drawline(97.00,18.00)(91.00,18.00)
\drawline(91.00,18.00)(91.00,03.00)
\drawline(91.00,03.00)(97.00,03.00)
\drawline(97.00,03.00)(97.00,-03.00)
\drawline(97.00,-03.00)(68.00,-03.00)
\drawline(68.00,08.00)(68.00,23.00)
\drawline(68.00,23.00)(93.00,32.00)
\drawline(93.00,32.00)(97.00,32.00)
\drawline(97.00,32.00)(97.00,18.00)
\drawline(68.00,28.00)(68.00,32.00)
\drawline(68.00,32.00)(74.00,32.00)
\drawline(74.00,32.00)(74.00,48.00)
\drawline(74.00,48.00)(68.00,48.00)
\drawline(68.00,48.00)(68.00,52.00)
\drawline(68.00,52.00)(97.00,52.00)
\drawline(97.00,52.00)(97.00,38.00)
\drawline(97.00,38.00)(93.00,38.00)
\drawline(93.00,38.00)(72.00,28.00)
\drawline(72.00,28.00)(68.00,28.00)
\drawline(68.00,68.00)(68.00,72.00)
\drawline(68.00,72.00)(97.00,72.00)
\drawline(97.00,72.00)(97.00,58.00)
\drawline(68.00,38.00)(68.00,42.00)
\drawline(68.00,42.00)(72.00,42.00)
\drawline(72.00,42.00)(72.00,38.00)
\drawline(72.00,38.00)(68.00,38.00)
\drawline(68.00,58.00)(68.00,68.00)
\drawline(97.00,58.00)(68.00,58.00)
\drawline(120.00,-03.00)(120.00,12.00)
\drawline(120.00,12.00)(124.00,12.00)
\drawline(124.00,12.00)(145.00,32.00)
\drawline(145.00,32.00)(149.00,32.00)
\drawline(149.00,32.00)(149.00,18.00)
\drawline(149.00,18.00)(143.00,18.00)
\drawline(143.00,18.00)(143.00,02.00)
\drawline(143.00,02.00)(149.00,02.00)
\drawline(149.00,02.00)(149.00,-03.00)
\drawline(149.00,-03.00)(120.00,-03.00)
\drawline(120.00,18.00)(120.00,42.00)
\drawline(120.00,42.00)(145.00,52.00)
\drawline(145.00,52.00)(149.00,52.00)
\drawline(149.00,52.00)(149.00,38.00)
\drawline(149.00,38.00)(145.00,38.00)
\drawline(145.00,38.00)(124.00,18.00)
\drawline(124.00,18.00)(120.00,18.00)
\drawline(120.00,48.00)(120.00,62.00)
\drawline(120.00,62.00)(145.00,72.00)
\drawline(145.00,72.00)(149.00,72.00)
\drawline(149.00,72.00)(149.00,58.00)
\drawline(149.00,58.00)(124.00,48.00)
\drawline(124.00,48.00)(120.00,48.00)
\drawline(120.00,68.00)(120.00,72.00)
\drawline(120.00,72.00)(124.00,72.00)
\drawline(124.00,72.00)(124.00,68.00)
\drawline(124.00,68.00)(120.00,68.00)
\drawline(145.00,08.00)(145.00,12.00)
\drawline(145.00,12.00)(149.00,12.00)
\drawline(149.00,12.00)(149.00,08.00)
\drawline(149.00,08.00)(145.00,08.00)
\drawline(93.00,08.00)(93.00,12.00)
\drawline(93.00,12.00)(97.00,12.00)
\drawline(97.00,12.00)(97.00,08.00)
\drawline(97.00,08.00)(93.00,08.00)
\put(70.00,00.00){\makebox(0,0)[cc]{$\bullet$}}
\put(70.00,10.00){\makebox(0,0)[cc]{$\bullet$}}
\put(70.00,20.00){\makebox(0,0)[cc]{$\bullet$}}
\put(70.00,30.00){\makebox(0,0)[cc]{$\bullet$}}
\put(70.00,40.00){\makebox(0,0)[cc]{$\bullet$}}
\put(70.00,50.00){\makebox(0,0)[cc]{$\bullet$}}
\put(70.00,60.00){\makebox(0,0)[cc]{$\bullet$}}
\put(70.00,70.00){\makebox(0,0)[cc]{$\bullet$}}
\put(95.00,10.00){\makebox(0,0)[cc]{$\bullet$}}
\put(95.00,20.00){\makebox(0,0)[cc]{$\bullet$}}
\put(95.00,30.00){\makebox(0,0)[cc]{$\bullet$}}
\put(95.00,40.00){\makebox(0,0)[cc]{$\bullet$}}
\put(95.00,50.00){\makebox(0,0)[cc]{$\bullet$}}
\put(95.00,60.00){\makebox(0,0)[cc]{$\bullet$}}
\put(95.00,70.00){\makebox(0,0)[cc]{$\bullet$}}
\put(95.00,00.00){\makebox(0,0)[cc]{$\bullet$}}
\put(122.00,00.00){\makebox(0,0)[cc]{$\bullet$}}
\put(122.00,10.00){\makebox(0,0)[cc]{$\bullet$}}
\put(122.00,20.00){\makebox(0,0)[cc]{$\bullet$}}
\put(122.00,30.00){\makebox(0,0)[cc]{$\bullet$}}
\put(122.00,40.00){\makebox(0,0)[cc]{$\bullet$}}
\put(122.00,50.00){\makebox(0,0)[cc]{$\bullet$}}
\put(122.00,60.00){\makebox(0,0)[cc]{$\bullet$}}
\put(122.00,70.00){\makebox(0,0)[cc]{$\bullet$}}
\put(147.00,10.00){\makebox(0,0)[cc]{$\bullet$}}
\put(147.00,20.00){\makebox(0,0)[cc]{$\bullet$}}
\put(147.00,30.00){\makebox(0,0)[cc]{$\bullet$}}
\put(147.00,40.00){\makebox(0,0)[cc]{$\bullet$}}
\put(147.00,50.00){\makebox(0,0)[cc]{$\bullet$}}
\put(147.00,60.00){\makebox(0,0)[cc]{$\bullet$}}
\put(147.00,70.00){\makebox(0,0)[cc]{$\bullet$}}
\put(147.00,00.00){\makebox(0,0)[cc]{$\bullet$}}
\drawline(45.00,00.30)(50.00,00.30)
\drawline(55.00,00.30)(60.00,00.30)
\drawline(65.00,00.30)(70.00,00.30)
\drawline(45.00,10.30)(50.00,10.30)
\drawline(55.00,10.30)(60.00,10.30)
\drawline(65.00,10.30)(70.00,10.30)
\drawline(45.00,20.30)(50.00,20.30)
\drawline(55.00,20.30)(60.00,20.30)
\drawline(65.00,20.30)(70.00,20.30)
\drawline(45.00,30.30)(50.00,30.30)
\drawline(55.00,30.30)(60.00,30.30)
\drawline(65.00,30.30)(70.00,30.30)
\drawline(45.00,40.30)(50.00,40.30)
\drawline(55.00,40.30)(60.00,40.30)
\drawline(65.00,40.30)(70.00,40.30)
\drawline(45.00,50.30)(50.00,50.30)
\drawline(55.00,50.30)(60.00,50.30)
\drawline(65.00,50.30)(70.00,50.30)
\drawline(45.00,60.30)(50.00,60.30)
\drawline(55.00,60.30)(60.00,60.30)
\drawline(65.00,60.30)(70.00,60.30)
\drawline(45.00,70.30)(50.00,70.30)
\drawline(55.00,70.30)(60.00,70.30)
\drawline(65.00,70.30)(70.00,70.30)
\put(57.50,35.00){\makebox(0,0)[cc]{$\cdot$}}
\put(109.50,35.00){\makebox(0,0)[cc]{$=$}}
\end{picture}
\caption{An example of a non-trivial product in 
$\widehat{\mathcal{PI}}^*_8$.}\label{figthree}
\end{figure}

There exist at least two different ways to deform the multiplication on the
semigroup  $\mathcal{PI}^*_k$. As we will now work with different 
multiplications, to distinguish them we denote by $\cdot$ the usual 
multiplication in $\mathcal{PI}^*_k$. The first ``naive'' deformation can be
constructed for any inverse semigroups (see e.g. \cite[4.1]{St}). For the 
semigroup $\mathcal{PI}^*_k$ this works as follows: Set 
$\widehat{\mathcal{PI}}^*_k=\mathcal{PI}^*_k\cup\{\mathbf{0}\}$. For
$\alpha,\beta\in \mathcal{PI}^*_k$ consider the following condition:
\begin{equation}\label{cond1}
\begin{array}{l}
\text{ If $A$ is a block of $\alpha$ and $B$ is a block of $\beta$ such that}\\
\text{ $(A\cap \mathbf{K}')\cap (B\cap \mathbf{K})'\neq \varnothing$,
then $A\cap \mathbf{K}'= (B\cap \mathbf{K})'$.}
\end{array}
\end{equation}
Define a new operation $\star$ on $\widehat{\mathcal{PI}}^*_k$ as follows:
\begin{displaymath}
\alpha\star\beta=
\begin{cases}
\alpha\cdot\beta,& \alpha,\beta\in \mathcal{PI}^*_k 
\text{ and \eqref{cond1} is satisfied},\\
\mathbf{0}, & \text{otherwise}.
\end{cases}
\end{displaymath}
For example, the $\star$-product of the two elements on the left hand side of
Figure~\ref{figtwo} equals $\mathbf{0}$. An example of a non-trivial
product in $\widehat{\mathcal{PI}}^*_k$ is shown on Figure~\ref{figthree}.

From the general theory (see e.g. \cite[4.1]{St}) it follows that 
$\widehat{\mathcal{PI}}^*_k$ is a semigroup. By \cite[Lemma~4.1]{St} 
the map
\begin{displaymath}
\begin{array}{rcl}
\varphi:\mathbb{C}[\mathcal{PI}^*_k]&\longrightarrow&
\overline{\mathbb{C}[\widehat{\mathcal{PI}}^*_k]},\\
\alpha&\longmapsto& \sum_{\beta\succeq \alpha}\beta
\end{array}
\end{displaymath}
is an algebra isomorphism. The map $\varphi$ allows one to reformulate
Theorem~\ref{thm2} in terms of the right action on $U^{\otimes k}$ of 
the semigroup $\widehat{\mathcal{PI}}^*_k$. This is fairly straightforward. 
What we would like to do is to give an explicit combinatorial description of 
the action of $\widehat{\mathcal{PI}}^*_k$ on $U^{\otimes k}$, which is 
induced by the action of $\mathcal{PI}^*_k$.  For  $\alpha\in
\widehat{\mathcal{PI}}^*_k\setminus\{\mathbf{0}\}$ and 
$\mathbf{i}\in (\mathbf{N}\cup\{0\})^k$ let 
$\hat{M}(\alpha,\mathbf{i})$ denote the set of all 
$\mathbf{l}\in (\mathbf{N}\cup\{0\})^k$ such that 
\begin{itemize}
\item for each block $\{a_1,\dots,a_p,b'_1,\dots,b'_q\}$ of $\alpha$ we have 
\begin{displaymath}
i_{a_1}=i_{a_2}=\dots=i_{a_p}=l_{b_1}=l_{b_2}=\dots=l_{b_q}\neq 0
\end{displaymath}
(we will say that this common value is the number, corresponding to this block);
\item numbers, which correspond to different blocks of $\alpha$, are different;
\item for any $a\in \mathbf{K}$ and $b\in \mathbf{K}'$ which do not belong 
to any block of $\alpha$ we have $i_a=0=l_b$.
\end{itemize}
If $\alpha=\mathbf{0}$, we set $\hat{M}(\alpha,\mathbf{i})=\varnothing$. As
before, it is easy to see that $|\hat{M}(\alpha,\mathbf{i})|\leq 1$ for all
$\alpha$ and $\mathbf{i}$. Set 
\begin{equation}\label{form4}
\alpha\star v_{\mathbf{i}}=
\begin{cases}
\sum_{\mathbf{l}\in \hat{M}(\alpha,\mathbf{i})} v_{\mathbf{l}},&
\hat{M}(\alpha,\mathbf{i})\neq \varnothing;\\
0, & \text{otherwise.} 
\end{cases}
\end{equation}

\begin{proposition}\label{prop3}
\begin{enumerate}[(i)]
\item\label{prop3.1} Let $\alpha\in {\mathcal{PI}}^*_k$ and $\mathbf{i}\in
(\mathbf{N}\cup\{0\})^k$. Then $\alpha(v_{\mathbf{i}})= 0$ implies
$\beta\star v_{\mathbf{i}}= 0$ for all $\beta\succeq\alpha$.
\item\label{prop3.2} Let $\alpha\in {\mathcal{PI}}^*_k$ and $\mathbf{i}\in
(\mathbf{N}\cup\{0\})^k$. Then $\alpha(v_{\mathbf{i}})\neq 0$ implies
that there exists a unique $\beta\succeq\alpha$ such that 
$\beta\star v_{\mathbf{i}}\neq  0$.
\item\label{prop3.3}
For any $\alpha\in \widehat{\mathcal{PI}}^*_k$ we have 
$\alpha\star v_{\mathbf{i}}=\varphi^{-1}(\alpha)(v_{\mathbf{i}})$.
\end{enumerate}
\end{proposition}

\begin{proof}
By \eqref{form1}, the condition $\alpha(v_{\mathbf{i}})= 0$ means that there 
exists some block $A$ of $\alpha$ and $a,b\in A$ such that $i_a\neq i_b$. If 
$\beta\succeq\alpha$, then there exists a block $B$ of $\beta$, which contains
$A$. We still have $i_a\neq i_b$ for $a,b\in B$. Hence 
$\beta\star v_{\mathbf{i}}= 0$ by \eqref{form4}. This proves \eqref{prop3.1}.

By \eqref{form1}, the condition $\alpha(v_{\mathbf{i}})\neq 0$ means that $i_a= 
i_b$ for all $a,b$ from the same block of $\alpha$. There is a unique way
to unite blocks of $A$ into bigger blocks such that the latter property still
holds for these bigger blocks, but $i_a\neq i_b$ if $a$ and $b$ are form 
different bigger blocks. By \eqref{form4}, this defines a unique 
$\beta\succeq\alpha$ such that  $\beta\star v_{\mathbf{i}}\neq  0$.

The statement \eqref{prop3.3} follows from \eqref{prop3.1}, \eqref{prop3.2}
and the definition of $\varphi$.
\end{proof}

We reamark that an explicit formula for $\varphi^{-1}(\alpha)$ can be obtained
using the M{\"o}bius inversion formula with respect to the partial order
$\preceq$ on the set ${\mathcal{PI}}^*_k$.

\begin{figure}
\special{em:linewidth 0.4pt} \unitlength 0.80mm
\linethickness{0.4pt}
\begin{picture}(150.00,75.00)
\put(20.00,00.00){\makebox(0,0)[cc]{$\bullet$}}
\put(20.00,10.00){\makebox(0,0)[cc]{$\bullet$}}
\put(20.00,20.00){\makebox(0,0)[cc]{$\bullet$}}
\put(20.00,30.00){\makebox(0,0)[cc]{$\bullet$}}
\put(20.00,40.00){\makebox(0,0)[cc]{$\bullet$}}
\put(20.00,50.00){\makebox(0,0)[cc]{$\bullet$}}
\put(20.00,60.00){\makebox(0,0)[cc]{$\bullet$}}
\put(20.00,70.00){\makebox(0,0)[cc]{$\bullet$}}
\put(45.00,10.00){\makebox(0,0)[cc]{$\bullet$}}
\put(45.00,20.00){\makebox(0,0)[cc]{$\bullet$}}
\put(45.00,30.00){\makebox(0,0)[cc]{$\bullet$}}
\put(45.00,40.00){\makebox(0,0)[cc]{$\bullet$}}
\put(45.00,50.00){\makebox(0,0)[cc]{$\bullet$}}
\put(45.00,60.00){\makebox(0,0)[cc]{$\bullet$}}
\put(45.00,70.00){\makebox(0,0)[cc]{$\bullet$}}
\put(45.00,00.00){\makebox(0,0)[cc]{$\bullet$}}
\drawline(93.00,08.00)(93.00,12.00)
\drawline(93.00,12.00)(97.00,12.00)
\drawline(97.00,12.00)(97.00,08.00)
\drawline(97.00,08.00)(93.00,08.00)
\drawline(18.00,-03.00)(18.00,13.00)
\drawline(18.00,13.00)(22.00,13.00)
\drawline(22.00,13.00)(43.00,23.00)
\drawline(43.00,23.00)(47.00,23.00)
\drawline(47.00,23.00)(47.00,-03.00)
\drawline(47.00,-03.00)(18.00,-03.00)
\drawline(18.00,18.00)(18.00,32.00)
\drawline(18.00,32.00)(47.00,32.00)
\drawline(47.00,32.00)(47.00,27.00)
\drawline(47.00,27.00)(22.00,18.00)
\drawline(22.00,18.00)(18.00,18.00)
\drawline(18.00,38.00)(18.00,42.00)
\drawline(18.00,42.00)(43.00,52.00)
\drawline(43.00,52.00)(47.00,52.00)
\drawline(47.00,52.00)(47.00,48.00)
\drawline(47.00,48.00)(22.00,38.00)
\drawline(22.00,38.00)(18.00,38.00)
\drawline(18.00,48.00)(18.00,62.00)
\drawline(18.00,62.00)(43.00,72.00)
\drawline(43.00,72.00)(47.00,72.00)
\drawline(47.00,72.00)(47.00,58.00)
\drawline(47.00,58.00)(22.00,48.00)
\drawline(22.00,48.00)(18.00,48.00)
\drawline(18.00,68.00)(18.00,72.00)
\drawline(18.00,72.00)(22.00,72.00)
\drawline(22.00,72.00)(22.00,68.00)
\drawline(22.00,68.00)(18.00,68.00)
\drawline(43.00,38.00)(43.00,42.00)
\drawline(43.00,42.00)(47.00,42.00)
\drawline(47.00,42.00)(47.00,38.00)
\drawline(47.00,38.00)(43.00,38.00)
\drawline(68.00,-03.00)(68.00,03.00)
\drawline(68.00,03.00)(91.00,22.00)
\drawline(91.00,22.00)(97.00,22.00)
\drawline(97.00,22.00)(97.00,18.00)
\drawline(97.00,18.00)(91.00,18.00)
\drawline(91.00,18.00)(91.00,03.00)
\drawline(91.00,03.00)(97.00,03.00)
\drawline(97.00,03.00)(97.00,-03.00)
\drawline(97.00,-03.00)(68.00,-03.00)
\drawline(68.00,18.00)(68.00,23.00)
\drawline(68.00,23.00)(93.00,32.00)
\drawline(93.00,32.00)(97.00,32.00)
\drawline(97.00,32.00)(97.00,28.00)
\drawline(97.00,28.00)(93.00,28.00)
\drawline(93.00,28.00)(72.00,18.00)
\drawline(72.00,18.00)(68.00,18.00)
\drawline(68.00,28.00)(68.00,32.00)
\drawline(68.00,32.00)(74.00,32.00)
\drawline(74.00,32.00)(74.00,48.00)
\drawline(74.00,48.00)(68.00,48.00)
\drawline(68.00,48.00)(68.00,52.00)
\drawline(68.00,52.00)(97.00,52.00)
\drawline(97.00,52.00)(97.00,38.00)
\drawline(97.00,38.00)(93.00,38.00)
\drawline(93.00,38.00)(72.00,28.00)
\drawline(72.00,28.00)(68.00,28.00)
\drawline(68.00,58.00)(68.00,72.00)
\drawline(68.00,72.00)(97.00,72.00)
\drawline(97.00,72.00)(97.00,68.00)
\drawline(97.00,68.00)(72.00,58.00)
\drawline(72.00,58.00)(68.00,58.00)
\drawline(68.00,38.00)(68.00,42.00)
\drawline(68.00,42.00)(72.00,42.00)
\drawline(72.00,42.00)(72.00,38.00)
\drawline(72.00,38.00)(68.00,38.00)
\drawline(93.00,58.00)(93.00,62.00)
\drawline(93.00,62.00)(97.00,62.00)
\drawline(97.00,62.00)(97.00,58.00)
\drawline(97.00,58.00)(93.00,58.00)
\drawline(68.00,08.00)(68.00,12.00)
\drawline(68.00,12.00)(72.00,12.00)
\drawline(72.00,12.00)(72.00,08.00)
\drawline(72.00,08.00)(68.00,08.00)
\drawline(120.00,48.00)(120.00,62.00)
\drawline(120.00,62.00)(145.00,72.00)
\drawline(145.00,72.00)(149.00,72.00)
\drawline(149.00,72.00)(149.00,68.00)
\drawline(149.00,68.00)(124.00,48.00)
\drawline(124.00,48.00)(120.00,48.00)
\drawline(120.00,68.00)(120.00,72.00)
\drawline(120.00,72.00)(124.00,72.00)
\drawline(124.00,72.00)(124.00,68.00)
\drawline(124.00,68.00)(120.00,68.00)
\drawline(120.00,38.00)(120.00,42.00)
\drawline(120.00,42.00)(124.00,42.00)
\drawline(124.00,42.00)(124.00,38.00)
\drawline(124.00,38.00)(120.00,38.00)
\drawline(120.00,28.00)(120.00,32.00)
\drawline(120.00,32.00)(124.00,32.00)
\drawline(124.00,32.00)(124.00,28.00)
\drawline(124.00,28.00)(120.00,28.00)
\drawline(120.00,18.00)(120.00,22.00)
\drawline(120.00,22.00)(124.00,22.00)
\drawline(124.00,22.00)(124.00,18.00)
\drawline(124.00,18.00)(120.00,18.00)
\drawline(120.00,08.00)(120.00,12.00)
\drawline(120.00,12.00)(124.00,12.00)
\drawline(124.00,12.00)(124.00,08.00)
\drawline(124.00,08.00)(120.00,08.00)
\drawline(120.00,-02.00)(120.00,02.00)
\drawline(120.00,02.00)(124.00,02.00)
\drawline(124.00,02.00)(124.00,-02.00)
\drawline(124.00,-02.00)(120.00,-02.00)
\drawline(145.00,08.00)(145.00,12.00)
\drawline(145.00,12.00)(149.00,12.00)
\drawline(149.00,12.00)(149.00,08.00)
\drawline(149.00,08.00)(145.00,08.00)
\drawline(145.00,58.00)(145.00,62.00)
\drawline(145.00,62.00)(149.00,62.00)
\drawline(149.00,62.00)(149.00,58.00)
\drawline(149.00,58.00)(145.00,58.00)
\drawline(145.00,-02.00)(145.00,02.00)
\drawline(145.00,02.00)(149.00,02.00)
\drawline(149.00,02.00)(149.00,-02.00)
\drawline(149.00,-02.00)(145.00,-02.00)
\drawline(145.00,18.00)(145.00,22.00)
\drawline(145.00,22.00)(149.00,22.00)
\drawline(149.00,22.00)(149.00,18.00)
\drawline(149.00,18.00)(145.00,18.00)
\drawline(145.00,28.00)(145.00,32.00)
\drawline(145.00,32.00)(149.00,32.00)
\drawline(149.00,32.00)(149.00,28.00)
\drawline(149.00,28.00)(145.00,28.00)
\drawline(145.00,38.00)(145.00,42.00)
\drawline(145.00,42.00)(149.00,42.00)
\drawline(149.00,42.00)(149.00,38.00)
\drawline(149.00,38.00)(145.00,38.00)
\drawline(145.00,48.00)(145.00,52.00)
\drawline(145.00,52.00)(149.00,52.00)
\drawline(149.00,52.00)(149.00,48.00)
\drawline(149.00,48.00)(145.00,48.00)
\put(70.00,00.00){\makebox(0,0)[cc]{$\bullet$}}
\put(70.00,10.00){\makebox(0,0)[cc]{$\bullet$}}
\put(70.00,20.00){\makebox(0,0)[cc]{$\bullet$}}
\put(70.00,30.00){\makebox(0,0)[cc]{$\bullet$}}
\put(70.00,40.00){\makebox(0,0)[cc]{$\bullet$}}
\put(70.00,50.00){\makebox(0,0)[cc]{$\bullet$}}
\put(70.00,60.00){\makebox(0,0)[cc]{$\bullet$}}
\put(70.00,70.00){\makebox(0,0)[cc]{$\bullet$}}
\put(95.00,10.00){\makebox(0,0)[cc]{$\bullet$}}
\put(95.00,20.00){\makebox(0,0)[cc]{$\bullet$}}
\put(95.00,30.00){\makebox(0,0)[cc]{$\bullet$}}
\put(95.00,40.00){\makebox(0,0)[cc]{$\bullet$}}
\put(95.00,50.00){\makebox(0,0)[cc]{$\bullet$}}
\put(95.00,60.00){\makebox(0,0)[cc]{$\bullet$}}
\put(95.00,70.00){\makebox(0,0)[cc]{$\bullet$}}
\put(95.00,00.00){\makebox(0,0)[cc]{$\bullet$}}
\put(122.00,00.00){\makebox(0,0)[cc]{$\bullet$}}
\put(122.00,10.00){\makebox(0,0)[cc]{$\bullet$}}
\put(122.00,20.00){\makebox(0,0)[cc]{$\bullet$}}
\put(122.00,30.00){\makebox(0,0)[cc]{$\bullet$}}
\put(122.00,40.00){\makebox(0,0)[cc]{$\bullet$}}
\put(122.00,50.00){\makebox(0,0)[cc]{$\bullet$}}
\put(122.00,60.00){\makebox(0,0)[cc]{$\bullet$}}
\put(122.00,70.00){\makebox(0,0)[cc]{$\bullet$}}
\put(147.00,10.00){\makebox(0,0)[cc]{$\bullet$}}
\put(147.00,20.00){\makebox(0,0)[cc]{$\bullet$}}
\put(147.00,30.00){\makebox(0,0)[cc]{$\bullet$}}
\put(147.00,40.00){\makebox(0,0)[cc]{$\bullet$}}
\put(147.00,50.00){\makebox(0,0)[cc]{$\bullet$}}
\put(147.00,60.00){\makebox(0,0)[cc]{$\bullet$}}
\put(147.00,70.00){\makebox(0,0)[cc]{$\bullet$}}
\put(147.00,00.00){\makebox(0,0)[cc]{$\bullet$}}
\drawline(45.00,00.30)(50.00,00.30)
\drawline(55.00,00.30)(60.00,00.30)
\drawline(65.00,00.30)(70.00,00.30)
\drawline(45.00,10.30)(50.00,10.30)
\drawline(55.00,10.30)(60.00,10.30)
\drawline(65.00,10.30)(70.00,10.30)
\drawline(45.00,20.30)(50.00,20.30)
\drawline(55.00,20.30)(60.00,20.30)
\drawline(65.00,20.30)(70.00,20.30)
\drawline(45.00,30.30)(50.00,30.30)
\drawline(55.00,30.30)(60.00,30.30)
\drawline(65.00,30.30)(70.00,30.30)
\drawline(45.00,40.30)(50.00,40.30)
\drawline(55.00,40.30)(60.00,40.30)
\drawline(65.00,40.30)(70.00,40.30)
\drawline(45.00,50.30)(50.00,50.30)
\drawline(55.00,50.30)(60.00,50.30)
\drawline(65.00,50.30)(70.00,50.30)
\drawline(45.00,60.30)(50.00,60.30)
\drawline(55.00,60.30)(60.00,60.30)
\drawline(65.00,60.30)(70.00,60.30)
\drawline(45.00,70.30)(50.00,70.30)
\drawline(55.00,70.30)(60.00,70.30)
\drawline(65.00,70.30)(70.00,70.30)
\put(57.50,35.00){\makebox(0,0)[cc]{$\cdot$}}
\put(109.50,35.00){\makebox(0,0)[cc]{$=$}}
\end{picture}
\caption{Elements of $\widetilde{\mathcal{PI}}^*_8$ and their
multiplication.}\label{figfour}
\end{figure}
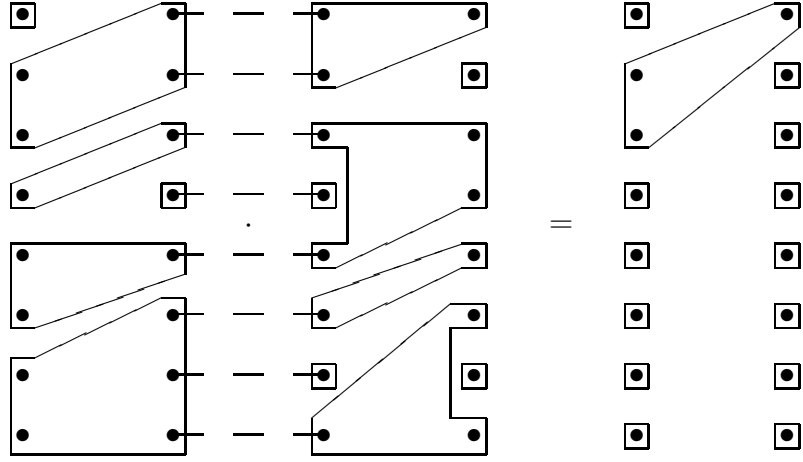

Another deformation of $\mathcal{PI}^*_k$ was proposed in \cite{Ve} and
studied in \cite{KMl}. The idea is rather similar to the ``naive'' deformation
$\widehat{\mathcal{PI}}^*_k$, however, the deformation of the multiplication 
in this second case is more subtle than in the ``naive'' case. Set
$\widetilde{\mathcal{PI}}^*_k=\mathcal{PI}^*_k$ (as a set). For 
$\alpha,\beta\in \widetilde{\mathcal{PI}}^*_k$ define the
element $\alpha\bullet\beta\in \widetilde{\mathcal{PI}}^*_k$ in the following 
way: A block $C$ belongs to $\alpha\bullet\beta$ if and only if there exists a 
block $A$ of $\alpha$ and a block $B$ of $\beta$ such that $A\cap
\mathbf{K}'=(B\cap \mathbf{K})'$ and $C=(A\cap \mathbf{K})\cup 
(B\cap \mathbf{K}')$. An example of multiplication for two elements from 
$\widetilde{\mathcal{PI}}^*_8$ is shown on Figure~\ref{figfour}.

For $\alpha,\beta\in \widetilde{\mathcal{PI}}^*_k$ we will write 
$\beta\vdash\alpha$ provided that each block of $\beta$ is a block of $\alpha$. 
By \cite[Lemma~4.1]{St} the mapping 
\begin{displaymath}
\begin{array}{rcl}
\psi:\mathbb{C}[\widetilde{\mathcal{PI}}^*_k]&\longrightarrow&
\overline{\mathbb{C}[\widehat{\mathcal{PI}}^*_k]},\\
\alpha&\longmapsto& \sum_{\beta\vdash\alpha}\beta
\end{array}
\end{displaymath}
is an algebra isomorphism. The maps $\psi$ and $\varphi$ allow one to reformulate Theorem~\ref{thm2} in terms of the right action on $U^{\otimes k}$
of the semigroup $\widetilde{\mathcal{PI}}^*_k$. This is again fairly
straightforward, so we just  give an explicit combinatorial description of 
the action of $\widetilde{\mathcal{PI}}^*_k$ on $U^{\otimes k}$, which is 
induced by the action of $\mathcal{PI}^*_k$.  For  $\alpha\in
\widetilde{\mathcal{PI}}^*_k$ and $\mathbf{i}\in
(\mathbf{N}\cup\{0\})^k$ let $\tilde{M}(\alpha,\mathbf{i})$ denote the set of 
all $\mathbf{l}\in (\mathbf{N}\cup\{0\})^k$ such that 
\begin{itemize}
\item for each block $\{a_1,\dots,a_p,b'_1,\dots,b'_q\}$ of $\alpha$ we have 
\begin{displaymath}
i_{a_1}=i_{a_2}=\dots=i_{a_p}=l_{b_1}=l_{b_2}=\dots=l_{b_q};
\end{displaymath}
\item non-zero numbers, which correspond to different blocks of $\alpha$, 
are different;
\item for any $a\in \mathbf{K}$ and $b\in \mathbf{K}'$ which do not belong 
to any block of $\alpha$ we have $i_a=0=l_b$.
\end{itemize}
Again we have $|\tilde{M}(\alpha,\mathbf{i})|\leq 1$ for all $\alpha$ 
and $\mathbf{i}$. Set 
\begin{displaymath}
\alpha\bullet v_{\mathbf{i}}=
\begin{cases}
\sum_{\mathbf{l}\in \tilde{M}(\alpha,\mathbf{i})} v_{\mathbf{l}},&
\tilde{M}(\alpha,\mathbf{i})\neq \varnothing;\\
0, & \text{otherwise.} 
\end{cases}
\end{displaymath}

\begin{proposition}\label{prop4}
For any $\alpha\in \widetilde{\mathcal{PI}}^*_k$ and  any
$\mathbf{i}\in (\mathbf{N}\cup\{0\})^k$ we have 
$\alpha\bullet v_{\mathbf{i}}=\psi(\alpha)\star v_{\mathbf{i}}$.
\end{proposition}

\begin{proof}
This follows immediately from the definitions.
\end{proof}

\vspace{0.5cm}

\noindent
V.M.: Department of Mathematics, Uppsala University, SE 471 06,
Uppsala, SWEDEN, e-mail: {\em mazor\symbol{64}math.uu.se}
\vspace{0.5cm}

\noindent
G.K.: Department of Mechanics and Mathematics, Kyiv Taras Shev\-chen\-ko
University, 64, Volodymyrska st., 01033, Kyiv, UKRAINE,
e-mail: {\em akudr\symbol{64}univ.kiev.ua}

\end{document}